\documentstyle[amscd,amssymb,verbatim,epsf]{amsart}
\newcommand{\R}{\mathbb{R}}
\newcommand{\N}{{\mathbb{N}}_{0}}

\begin{document}

\theoremstyle{plain}
\newtheorem{Thm}{Theorem}
\newtheorem{Cor}{Corollary}
\newtheorem{Con}{Conjecture}
\newtheorem{Main}{Main Theorem}

\newtheorem{Lem}{Lemma}
\newtheorem{Prop}{Proposition}

\theoremstyle{definition}
\newtheorem{Def}{Definition}
\newtheorem{Note}{Note}

\theoremstyle{acknowledgments}
\newtheorem{Acknow}{Acknowledgments}

\theoremstyle{remark}
\newtheorem{notation}{Notation}

\errorcontextlines=0
\numberwithin{equation}{section}

\title[On Abel statistical convergence]%
   {On Abel statistical convergence}
\author[Iffet Taylan]{I. Taylan and H. Cakalli\\
Maltepe University, \.{I}stanbul-Turkey \\Phone:(+90216)6261050,fax:(+90216)6261113}

\address{Iffet Taylan,  \; \; \; Maltepe University, Education Faculty, Marmara E\u{g}\.{I}t\.{I}m K\"oy\"u, TR 34857, Maltepe, \.{I}stanbul-Turkey Phone: (+90216)6261050 Ext:2230, FA\textbf{x}: (+90216)6261093
}
\email{iffettaylan@@maltepe.edu.tr}

\address{H\"usey\.{I}n \c{C}akall\i \\
          Maltepe University, Graduate School of Science and Engineering, Department of Mathematics, Marmara E\u{g}\.{I}t\.{I}m K\"oy\"u, TR 34857, Maltepe, \.{I}stanbul-Turkey. Phone:(+90216)6261050 ext:2311, fax:(+90216)6261113}
\email{hcakalli@@maltepe.edu.tr; hcakalli@@istanbul.edu.tr;  hcakalli@@gmail.com}

\date{\today}

\maketitle

\begin{abstract}

In this paper, we introduce and investigate a concept of Abel statistical continuity. A real valued function $f$ is Abel statistically continuous on a subset $E$  of $\R$, the set of real numbers, if it preserves Abel statistical convergent sequences, i.e.  $(f(p_{k}))$ is Abel statistically convergent whenever $(p_{k})$ is an Abel statistical convergent sequence of points in $E$, where a sequence $(p_{k})$ of point in $\R$ is called Abel statistically convergent to a real number $L$  if Abel density of the set $\{k\in{\N}: |p_{k}-L|\geq\varepsilon \}$ is $0$ for every $\varepsilon>0$. Some other types of continuities are also studied and interesting results are obtained.
\end{abstract}

\maketitle

\section{Introduction}

The concept of continuity and any concept involving continuity play a very important role not only in pure mathematics but also in other branches of sciences
involving mathematics especially in computer sciences, information theory, biological science, economics, and dynamical systems.

Throughout the paper, $\N$, and $\R$ will denote the set of non negative integers and the set of real numbers, respectively. The boldface letters $\bf{p}$, $\textbf{q}$, $\bf{r}$, and $\bf{w}$ will be used for sequences $\textbf{p}=(p_{k})$, $\textbf{q}=(q_{k})$, $\textbf{r}=(r_{k})$, $\textbf{w}=(w_{k})$, ... of points in $\R$. A function $f:\R \longrightarrow \R$ is continuous if and only if it preserves convergent sequences. Using the idea of continuity of a real function in this manner, many kinds of continuities were introduced and investigated, not all but some of them we recall in the following: slowly oscillating continuity (\cite{CakalliSlowlyoscillatingcontinuity}, \cite{Vallin}, $\Delta$ slowly oscillating continuity (\cite{CanakandDikNewtypesofcontinuities}), ward continuity (\cite{CakalliForwardcontinuity}), statistical continuity,  (\cite{CakalliStatisticalwardcontinuity}), $\lambda$-statistical continuity (\cite{CakalliandSonmezandArasLamdastatisticallywardcontinuity}), lacunary statistical continuity (\cite{CakalliandArasandSonmezLacunarystatisticalwardcontinuity}), ideal sequential continuity (\cite{CakalliandHazarikaIdealquasiCauchysequences, CakalliAvariationonwardcontinuity}), and Abel sequential continuity (\cite{CakalliandAlbayrakNewtypecontinuitiesviaAbelconvergence}).
A real sequence $(p_{k})$ of points in $\mathbb{R}$ is called statistically convergent to an $L\in{\R}$ if $\lim_{n \to \infty}\frac{1}{n}|\{k\leq n: |p_{k}-L|\geq \varepsilon\}|=0$ for every $\varepsilon>0$, and this is denoted by $st-\lim p_{k}=L$ (\cite{Fast}, \cite{Fridy}, \cite{CakalliAstudyonstatisticalconvergence}, \cite{CasertaandKocinacOnstatisticalexhaustiveness}, \cite{CakalliAvariationonstatisticalwardcontinuity}, \cite{CakalliandKhan}, and  \cite{CasertaMaioKocinacStatisticalConvergenceinFunctionSpaces}). A sequence $(p_{k})$ is called lacunary statistically convergent (\cite{FridyandOrhanlacunarystatisconvergence}) to an $L\in{\R}$ if $\lim_{r\rightarrow\infty}\frac{1}{h_{r}}|\{k\in I_{r}: |p_{k}-L| \geq{\varepsilon}\}|=0$ for every $\varepsilon >0$, where $I_{r}=(k_{r-1},k_{r}]$, and $k_{0}=0$, $h_{r}:k_{r}-k_{r-1}\rightarrow \infty$ as $r\rightarrow\infty$ and $\theta=(k_{r})$ is an increasing sequence of positive integers, and this is denoted by $S_{\theta}-\lim p_{n}=L$ (see also \cite{CakalliLacunarystatisticalconvergenceintopgroups}). Throughout this paper, we assume that $\liminf_{r}\frac{k_{r}}{k_{r-1}}>1$.

The purpose of this paper is to investigate the concept of Abel statistical continuity, and present interesting results.

\section{Abel statistical continuity}
Recently the concept of Abel statistical convergence of a sequence is introduced and investigated in \cite{UnverAbelsummabilityintopologicalspaces} (see also \cite{PowellandShahSummabilityTheoryandItsApplications}). Although the definitions and the most of the results are also valid in a topological Hausdorff group, which allows countable local base at the origin, we investigate the notion in the real setting: a sequence $\textbf{p}=(p_{k})$ is called Abel statistically convergent to a real number $L$ if Abel density of the set $\{k\in{\N}: |p_{k}-L|\geq\varepsilon \}$ is $0$, i.e. $\lim_{x \to 1^{-}}(1-x)\sum_{k\in{\N}:|p_{k}-L|\geq\varepsilon}^{}x^{k}=0$ \\for every $\varepsilon>0$, and denoted by $Abel_{st}-lim \;p_{k}=L$ (\cite{PowellandShahSummabilityTheoryandItsApplications}). We note that Abel statistical limit of an Abel statistical convergent sequence is unique, i.e.
if a sequence $\textbf{p}=(p_{k})$ is Abel statistically convergent to $L_{1}$ and $L_{2}$, then $L_{1}=L_{2}$. Any convergent sequence is Abel statistically convergent, i.e. the  Abel statistical method is regular. On the other hand, the sum of two Abel statistical convergent sequence is Abel statistically convergent, that is if $(p_{k})$ and $(q_{k})$ are Abel statistical convergent sequences with $Abel_{st}-lim\; p_{k}=L_{1}$ and $Abel_{st}-lim\; q_{k}=L_{2}$, then $Abel_{st}-lim\; (p_{k}+q_{k})=Abel_{st}-lim\; p_{k} + Abel_{st}-lim \;q_{k}=L_{1} + L_{2}$. For any constant $c\in{\R}$, the sequence $(c.p_{k})$ is Abel statistically convergent whenever $(p_{k})$ is. Thus the set of Abel statistical convergent sequences forms a linear space. We note that Abel statistical convergence coincides with statistical convergence.
Abel statistical sequential method is a subsequential method in the manner of  \cite{CakalliSequentialdefinitionsofcompactness}, \cite{CakalliOnGcontinuity}, \cite{CakalliSequentialdefinitionsofconnectedness}, and  \cite{MucukandSahanOnGsequentialcontinuity}.
\begin{Def} \label{DefAbeþstatisticalcompactness}
A subset $E$ of $\mathbb{R}$ is called Abel statistically compact if any sequence of points in $E$ has an Abel statistical convergent subsequence whose Abel statistical limit is in $E$, i.e. whenever $\textbf{p}=(p_{n})$ is a sequence of points in $E$, there is an Abel statistical convergent subsequence $\textbf{r}=(r_{k})=(r_{n_{k}})$ of the sequence $\textbf{p}$ satisfying $Abel_{st}-lim\;\textbf{r}\in{E}$.
\end{Def}
According to this definition, any closed and bounded subset of $\R$ is Abel statistically compact, i.e. compact subsets of $\R$  are Abel statistically compact.
Notice that the union of two Abel statistical compact subsets of $\mathbb{R}$ is Abel statistically compact, therefore we see inductively that any finite union of Abel statistical compact subsets of $\mathbb{R}$ is Abel statistically compact, whereas the union of the infinite family of Abel statistical compact subsets of $\mathbb{R}$ is not always Abel statistical compact. These observations above suggest to us the following.
\begin{Thm} \label{TheoAbelstatisticalsequentialcompactcoincidewithboundednessandclosednes}
A subset of $\mathbb{R}$ is Abel statistically compact if and only if it is compact in the ordinary sense.
\end{Thm}
\begin{pf}
The proof can easily obtained, so is omitted.
\end{pf}

We note that Abel statistical sequential compactness coincides with not only with ordinary compactness, but also statistical sequential compactness (\cite{CakalliStatisticalwardcontinuity}), $\lambda$-statistical sequential compactness (\cite{CakalliandSonmezandAraslamdastatisticallywardcontinuity}), lacunary statistical sequential compactness (\cite{CakalliandArasandSonmezOnlacunarystatisticallyquasiCauchysequences}), strongly lacunary sequential compactness (\cite{CakalliNthetawardcontinuity}), Abel sequential compactness (\cite{CakalliandAlbayrakNewtypecontinuitiesviaAbelconvergence}), and $I$-sequential compactness for a nontrivial admissible ideal $I$ (\cite{CakalliandHazarikaIdealquasiCauchysequences}, \cite{CakalliAvariationonwardcontinuity}).

\begin{Def}
A real number $L$ is said to be in the Abel statistical sequential closure of a subset $E$ of $\mathbb{R}$, denoted by $\overline{E}^{{Abel}_{st}}$ if there is a sequence $\textbf{p}=(p_{k})$ of points in $E$ such that $Abel_{st}-lim p_{k}=L$, and it is called Abel statistically sequentially closed if $\overline{E}^{{Abel}_{st}}=E$.
\end{Def}
Note that the preceding definition is a special case of the definition of $G$-sequential closure in \cite{CakalliSequentialdefinitionsofcompactness}. It is clear that $\overline{\phi}^{{Abel}_{st}}=\phi$ and $\overline{\mathbb{R}}^{{Abel}_{st}}=\mathbb{R}$.
It is easily seen that $E\subset \overline{E}\subset \overline{E}^{{Abel}_{st}}$.
We note that any Abel statistical sequentially closed subset of Abel statistical compact subset of $\mathbb{R}$ is also Abel statistical sequentially compact.
The condition closedness is essential here, in other words, a subset of an Abel statistical compact subset need not to be Abel statistically sequentially compact. For example the interval $]-1,1]$, i.e. the set of real numbers strictly greater than $-1$, and less than or equal to $1$, is a subset of Abel statistical compact subset $[-1,1]$, i.e. the set of real numbers greater than or equal to $-1$, and less than or equal to $1$, but not Abel statistically sequentially compact.
The intersection of two Abel statistically sequentially compact is Abel statistically sequentially compact. In general, the intersection of any family of Abel statistical sequentially compact subsets of $\R$ is Abel statistically sequentially compact.

We now introduce a new type of continuity defined via Abel statistical convergent sequences.
\begin{Def} \label{DefAbelstatisticalcontinuity}
A function $f$ is called Abel statistically continuous if it preserves Abel statistical convergent sequences, i.e. $(f(p_{k}))$ is Abel statistical convergent to  $f(L)$ whenever $(p_{k})$ is Abel statistically convergent to $L$.
\end{Def}
We see in the following that the sum of two Abel statistical continuous functions is Abel statistically continuous.
\begin{Prop}
If $f$ and $g$ are Abel statistical continuous functions on subset $E$ of $\R$, then $f+g$ is Abel statistically continuous on $E$.
\end{Prop}
\begin{pf}
The proof can be obtained easily so is omitted.
\end{pf}
If $f$ is Abel statistical continuous function, then $cf$ is Abel statistically continuous for any $c\in{\R}$. Thus the set of Abel statistical continuous functions forms a vector space. The composite of two Abel statistical continuous functions is Abel statistically continuous, and the product of two Abel statistical continuous functions is Abel statistically continuous. Abel statistical continuous image of $G$-sequentially connected  subset of $\R$ is $G$-sequentially connected (\cite{MucukandCakalliGsequentiallyconnectednessfortopologicalgroupswithoperations}). 

In connection with Abel statistical convergent sequences and convergent sequences the problem arises to investigate the following types of continuity of functions on $\mathbb{R}$.
\begin{description}
\item[($A_{st}$)] $(p_{n}) \in {\textbf{A}_{st}} \Rightarrow (f(p_{n})) \in {\textbf{A}_{st}}$
\item[$(A_{st}c)$] $(p_{n}) \in {\textbf{A}_{st}} \Rightarrow (f(p_{n})) \in {c}$
\item[$(c)$] $(p_{n}) \in {c} \Rightarrow (f(p_{n})) \in {c}$
\item[$(cA_{st})$] $(p_{n}) \in {c} \Rightarrow (f(p_{n})) \in {\textbf{A}_{st}}$
\end{description}
We see that $A_{st}$ is Abel statistical continuity of $f$, and $(c)$ states the ordinary continuity of $f$. We easily see that (c) implies (c$A_{st}$),  ($A_{st}$) implies (c$A_{st}$), and ($A_{st}$c) implies ($A_{st}$).  The converses are not always true as the identity function could be taken as a counter example for all the cases.

Now we give the implication $(A)$ implies $(c)$, i.e. any Abel statistical continuous function is continuous in the ordinary sense.
\begin{Thm}[] \label{TheoAbelstatisticalcontimpliescontinuity} If a function $f$ is Abel statistically continuous on a subset $E$ of $\mathbb{R}$, then it is continuous on $E$ in the ordinary sense.
\end{Thm}
\begin{pf}
Suppose that a function $f$ is not continuous on $E$. Then there exists a sequence $(p_{n})$ with $\lim_{n\to \infty } p_{n}=L$ such that $(f(p_{n}))$ is not
convergent to $f(L)$. If $(f(p_{n}))$ exists, and $\lim f(p_{n})$ is different from $f(L)$ then we easily see a contradiction.  Now suppose that $(f(p_{n}))$ has two subsequences of $f(p_{n})$ such that $\lim_{m \to \infty}f(p_{k_{m}})=L_{1}$ and  $\lim_{k \to \infty}f(p_{n_{k}})=L_{2}$. Since $(p_{n_{k}})$ is a subsequence of $(p_{n})$,  by hypothesis, $\lim_{k \to \infty}f(p_{n_{k}})=f(L_{1})$ and $(p_{k_{m}})$ is a subsequence of $(p_{n})$, by hypothesis  $\lim_{m \to \infty}f(p_{k_{m}})=f(L_{2})$. This leads to a contradiction.
If $(f(p_{n}))$ is unbounded above. Then we can find an $n_1$ such that $f(p_{n_{1}})>f(L)+2^{1}$. There exists a positive integer an $n_{2}>n_{1}$ such that $f(p_{n_{2}})>f(L)+2^{2}$. Suppose that we have chosen an $n_{k-1}>n_{k-2}$ such that $f(p_{n_{k-1}})>f(L)+2^{k-1}$. Then we can choose an $n_{k}>n_{k-1}$ such that $f(p_{n_{k}})>f(L)+2^{k}$. Inductively we can construct a subsequence $(f(p_{n_{k}}))$ of $(f(p_{n}))$ such that $f(p_{n_{k}})>f(L)+2^{k}$. Since the sequence $(p_{n_{k}})$ is a subsequence of  $(p_{n})$, the subsequence $(p_{n_{k}})$ is convergent so is Abel statistically convergent. But $(f(p_{n_{k}}))$ is not Abel statistically convergent as we see line below. For each positive integer $k$ we have $f(p_{n_{k}})>f(L)+2^{k}$. Thus \\
$\lim_{x \to 1^{-}}(1-x)\sum_{k\in{\N}:|f(p_{k})-f(L)|\geq \frac{1}{2}}^{}x^{k}=\lim_{x \to 1^{-}}(1-x)\sum_{k=0}^{\infty}x^{k}=\\\lim_{x \to 1^{-}}\frac{1-x}{1-x}=1\neq0$.
\end{pf}

We note that Abel statistical continuity implies not only ordinary continuity, but also statistical continuity, which follows from Theorem \ref{TheoAbelstatisticalcontimpliescontinuity}, \cite[Corollary 4]{CakalliAstudyonstatisticalconvergence}, Lemma 1 and Theorem 8 in \cite{CakalliOnGcontinuity},  lacunary statistical sequential continuity, which follows from Theorem \ref{TheoAbelstatisticalcontimpliescontinuity} (see also   \cite{CakalliandArasandSonmezOnlacunarystatisticallyquasiCauchysequences});  $\lambda$-statistical continuity, which follows from Theorem \ref{TheoAbelstatisticalcontimpliescontinuity}, Lemma 1 and Theorem 8 in \cite{CakalliOnGcontinuity} (see also \cite{CakalliandSonmezandAraslamdastatisticallywardcontinuity}); $N_{\theta}$-sequential continuity (\cite{CakalliNthetawardcontinuity}); $I$-sequential continuity for any non trivial admissible ideal (\cite[Theorem 4]{CakalliandHazarikaIdealquasiCauchysequences}); $G$-sequential continuity for any regular subsequential method $G$, which follows from Theorem 8 in  \cite{CakalliOnGcontinuity} (see also \cite{MucukandSahanOnGsequentialcontinuity})
\begin{Cor} Any Abel statistical continuous function on an Abel statistical compact subset of $\mathbb{R}$ is uniformly continuous.
\end{Cor}

It is well known that uniform limit of a sequence of continuous functions is continuous. This is also true for Abel statistical continuity, i.e. uniform limit of a sequence of Abel statistical continuous functions is Abel statistical continuous.
\begin{Thm} \label{Theorem 8}
If $(f_{n})$ is a sequence of Abel statistical continuous functions defined on a subset E of $\mathbb{R}$ and $(f_{n})$ is uniformly convergent to a function $f$, then $f$ is Abel statistically continuous on $E$.
\end{Thm}
\begin{pf}
Let $(p_{n})$ be an Abel statistical convergent sequence of real numbers in $E$. Write $Abel_{st}-\lim p_{n}=L.$ Take any $\varepsilon>0.$  Since $(f_{n})$ is uniformly convergent to $f$, there exists an  $N\in{\N}$ such
that
$ |f_{k}(t)-f(t)|<\varepsilon/3$ for all $t\in E$ whenever $k\geq N$.
Thus
$  \lim_{x \to 1^{-}}(1-x)\sum_{|f(p_{k})-f_{N}(p_{k})|\geq \frac{\varepsilon}{3}}^{}x^{k}=0$. Since $f_{N}$ is Abel statistically continuous, we have\\
$\lim_{x \to 1^{-}}(1-x) \sum_{|f_{N}(p_{k}))-f_{N}(L)|\geq \frac{\varepsilon}{3}}^{}x^{k}=0$.\\
On the other hand,\\
$ \sum_{|f(p_{k})-f(L)|\geq \varepsilon}^{}x^{k}\leq
\sum_{|f(p_{k})-f_{N}(p_{k})|\geq \frac{\varepsilon}{3}}^{}x^{k}+
\sum_{|f_{N}(p_{k}))-f_{N}(L)|\geq \frac{\varepsilon}{3}}^{}x^{k}+
\sum_{|f_{N}(L)-f(L)|\geq \frac{\varepsilon}{3}}^{}x^{k}$\\
for every $x$ satisfying $0<x<1$. Hence\\
$\lim_{x \to 1^{-}}(1-x) \sum_{|f(p_{k})-f(L)|\geq \varepsilon}^{}x^{k}\leq \\
\lim_{x \to 1^{-}}(1-x)\sum_{|f(p_{k})-f_{N}(p_{k})|\geq \frac{\varepsilon}{3}}^{}x^{k}+
\lim_{x \to 1^{-}}(1-x) \sum_{|f_{N}(p_{k}))-f_{N}(L)|\geq \frac{\varepsilon}{3}}^{}x^{k}+\\
\lim_{x \to 1^{-}}(1-x) \sum_{|f_{N}(L)-f(L)|\geq \frac{\varepsilon}{3}}^{}x^{k}=0+0+0=0$.\\
This completes the proof of the theorem.
\end{pf}
In the following theorem we prove that the set of Abel statistical continuous functions is a closed subset of the space of continuous functions.
\begin{Thm} \label{TheoThesetofAbelstatisticalcontinuousfunctionsisclosed}
The set of Abel statistical continuous functions on a subset $E$ of $\mathbb{R}$ is a closed subset of the set of all continuous functions on $E$, i.e.
$\overline{\textbf{A}_{st}C(E)}=\textbf{A}_{st}C(E)$, where $A_{st}C(E)$ is the set of all Abel statistical continuous functions on $E$, $\overline{\textbf{A}_{st}C(E)}$ denotes the set of all cluster points of $\textbf{A}_{st}C(E)$.
\end{Thm}
\begin{pf}
Let f be any element in the closure of $\textbf{AC}(E)$. Then there exists a sequence $(f_{n})$ of points in \textbf{AC}(E) such that $\lim_{k \to \infty}f_{k}=f.$  To show that f is Abel statistical continuous, take any Abel statistical convergent sequence  $(p_{k})$ of points $E$ with Abel statistical limit $L$. Let $\varepsilon>0$. Since $(f_{k})$ is convergent to $f$, there exists a positive integer $N$ such that $|f_{k}(t)-f(t)|<\varepsilon/3$ for all $t\in E$ whenever $k\geq N$. Thus
$  \lim_{x \to 1^{-}}(1-x)\sum_{|f(p_{k})-f_{N}(p_{k})|\geq \frac{\varepsilon}{3}}^{}x^{k}=0$. Since $f_{N}$ is Abel statistically continuous,
$\lim_{x \to 1^{-}}(1-x) \sum_{|f_{N}(p_{k}))-f_{N}(L)|\geq \frac{\varepsilon}{3}}^{}x^{k}=0$.
On the other hand,\\
$ \sum_{|f(p_{k})-f(L)|\geq \varepsilon}^{}x^{k}\leq \\
\sum_{|f(p_{k})-f_{N}(p_{k})|\geq \frac{\varepsilon}{3}}^{}x^{k} +
\sum_{|f_{N}(p_{k}))-f_{N}(L)|\geq \frac{\varepsilon}{3}}^{}x^{k} +
\sum_{|f_{N}(L)-f(L)|\geq \frac{\varepsilon}{3}}^{}x^{k}$\\
for every $x$ satisfying $0<x<1$. Hence\\
$\lim_{x \to 1^{-}}(1-x) \sum_{|f(p_{k})-f(L)|\geq \varepsilon}^{}x^{k}\leq \\
\lim_{x \to 1^{-}}(1-x)\sum_{|f(p_{k})-f_{N}(p_{k})|\geq \frac{\varepsilon}{3}}^{}x^{k} +
\lim_{x \to 1^{-}}(1-x) \sum_{|f_{N}(p_{k}))-f_{N}(L)|\geq \frac{\varepsilon}{3}}^{}x^{k} +\\
\lim_{x \to 1^{-}}(1-x) \sum_{|f_{N}(L)-f(L)|\geq \frac{\varepsilon}{3}}^{}x^{k}=0+0+0=0$. \\
This completes the proof of the theorem.
\end{pf}
\begin{Cor} \label{Corollary 10}
The set of all Abel statistical continuous functions on a subset $E$ of $\R$ is a complete subspace of the space of all continuous functions on $E$.
\end{Cor}
\begin{pf}
The proof follows from the preceding theorem, and the fact that the set of all continuous functions on $E$ is complete.
\end{pf}

\begin{Thm} \label{Theorem 11}
Abel statistical continuous image of any Abel statistical compact subset of $\R$ is Abel statistically sequentially compact.
\end{Thm}

\begin{pf}
Although the proof follows from \cite[Theorem 7]{CakalliSequentialdefinitionsofcompactness}, we give a short proof for completeness.
Let $f$ be any Abel statistical continuous function defined on a subset $E$ of $\mathbb{R}$ and $F$ be any Abel statistical compact subset of $E$. Take any sequence $\textbf{w}=(w_{k})$ of point in $f(F)$. Write $w_{k}=f(p_{k})$ for each $n\in{\N}$. Since $E$ is Abel statistically sequentially compact, there exists an Abel statistical convergent subsequence $\textbf{r}=(r_{k})$ of the sequence $\textbf{p}$. Write $Abel-lim\textbf{r}=L$. Since $f$ is Abel statistical continuous $Abel-limf(\textbf{r})=f(L)$. Thus $f(\textbf{r})=(f(r_{k}))$ is Abel statistically convergent to $f(L)$ which is a subsequence of the sequence $\textbf{w}$. This completes the proof.
\end{pf}
For $G:=Abel-lim$,  we have the following:

\begin{Thm} \label{TheoClosureofasetisasubestofclosureofimage}
If a function $f$ is Abel statistical continuous on a subset $E$ of $\R$, then $$f(\overline{B}^{Abel_{st}})\subset{\overline{(f(B))}^{Abel_{st}}}$$ for every subset $B$ of $E$.
\end{Thm}

\begin{pf}
The proof follows from the regularity of Abel  statistical sequential method, and Theorem 8 on page 316 of \cite{CakalliOnGcontinuity}.
\end{pf}

\maketitle
\section{Conclusion}

In this paper we introduce a concept of Abel statistical continuity, and present theorems related to this kind of continuity, and some other kinds of continuities. The concept of Abel statistical compactness is also introduced and investigated. One may expect this investigation to be a useful tool in the field of analysis in modeling various problems occurring in many areas of science, dynamical systems, computer science, information theory, and biological science. On the other hand, we suggest to introduce a concept of fuzzy Abel statistical convergence of sequences of fuzzy points, and investigate fuzzy Abel statistical continuity for fuzzy functions (see \cite{CakalliandPratulFuzzycompactnessviasummability}, \cite{KocinacSelectionpropertiesinfuzzymetricspaces} for the definitions and  related concepts in fuzzy setting). However due to the change in settings, the definitions and methods of proofs will not always be the same. An investigation of Abel sequential continuity and Abel sequential compactness can be done for double sequences (see \cite{CakalliandSavasStatisticalconvergenceofdoublesequences}, \cite{DjurcicandKocinacandZizovicDoublesequencesandselections}, \cite{CakalliandPattersonFunctionspreservingslowlyoscillatingdoublesequences} for basic concepts in the double sequences case). For some further study, we suggest to investigate Abel statistical quasi Cauchy sequences of points in a topological vector space valued cone metric space (see \cite{CakalliandSonmezandGenc}, \cite{PalandSavasandCakalliIconvergenceonconemetricspaces},  \cite{CakalliandSonmezSlowlyoscillatingcontinuityinabstractmetricspaces}, \cite{SonmezOnparacompactnessinconemetricspaces}, and \cite{SonmezandCakalliConenormedspacesandweightedmeans}) or in 2-normed spaces (\cite{MursaleenandAlotaibiOnIconvergenceinrandom2normedspaces}, \cite{CakalliandErsanStronglylacunarywardcontinuityin2normedspaces}, \cite{CakalliandErsanLacunarywardcontinuityin2normedspaces}, \cite{ErsanandCakalliWardContinuityin2NormedSpaces}).


\begin{thebibliography}{99}
\bibitem{Abel}  N.H. Abel,  Recherches sur la série 1+$\frac{m}{1}x$+$\frac{m(m-1)}{1.2}x^{2}$+.... J. Reine Angew. Math. \textbf{1} (1826) 311-339.
\bibitem{BadiozzamanandThorpeSomebest}   A.J. Badiozzaman, and B. Thorpe, Some best possible Tauberian results for Abel and Cesaro summability. Bull. London Math. Soc. \textbf{28} 3 (1996) 283-290.
\bibitem{BurtonandColemanQuasiCauchysequences} D. Burton, and J. Coleman,  Quasi-Cauchy Sequences, Amer. Math. Monthly. \textbf{117} 4 (2010) 328-333.
\bibitem{CakalliNthetawardcontinuity} H. Cakalli,  N-theta-Ward continuity, Abstr. Appl. Anal. \textbf{2012}, Article ID 680456 (2012) 8 pages   doi:10.1155/2012/680456.
\bibitem{CakalliAvariationonstatisticalwardcontinuity} H. Cakalli, A variation on statistical ward continuity, Bull. Malays. Math. Sci. Soc.  \textbf{\textbf{40}} (2017) 1701-1710. https://doi.org/10.1007/s40840-015-0195-0
\bibitem{CakalliVariationsonstatisticalquasiCauchysequences} H. Cakalli, Variations on statistical quasi Cauchy sequences, Bol. Soc. Paran. Mat. to appear.
\bibitem{CakalliandAlbayrakNewtypecontinuitiesviaAbelconvergence} H. Cakalli,  and M. Albayrak, New type continuities via Abel convergence. Scientific World Journal. \textbf{2014}, Article ID 398379  6 pages (2014) http://dx.doi.org/10.1155/2014/398379
\bibitem{CakalliandArasandSonmezLacunarystatisticalwardcontinuity} H. Cakalli, C.G. Aras, A. Sonmez, Lacunary statistical ward continuity, AIP Conf. Proc. \textbf{1676}, 020042 (2015); http://dx.doi.org/10.1063/1.4930468 arXiv \textbf{2013}:1102.1531
\bibitem{CakalliandHazarikaIdealquasiCauchysequences} H. \c{C}akall\i,  and B. Hazarika, Ideal quasi-Cauchy sequences, J. Inequal. Appl. \textbf{2012} (2012) 2012:234 ; doi:10.1186/1029-242X-2012-234
\bibitem{CakalliandPattersonFunctionspreservingslowlyoscillatingdoublesequences} H. Cakalli,   R.F. Patterson, Functions preserving slowly oscillating double sequences, An. Stiint. Univ. Al. I. Cuza Iasi. Mat. (N.S.) Tomul LXII, f. 2, vol. 2  (2016) 531-536. http://www.math.uaic.ro/~annalsmath/pdf-uri\%20anale/F2-2(2016)/Cakalli-Patterson.pdf
\bibitem{CakalliandTaylanAvariationonAbelstatisticalwardcontinuity} H. Cakalli and I. Taylan, A variation on Abel statistical ward continuity, AIP Conf. Proc.  \textbf{1676}, 020076 (2015) http://dx.doi.org/10.1063/1.4930502
\bibitem{CasertaandKocinacOnstatisticalexhaustiveness}  A. Caserta and Ljubisa D.R. Kocinac,  On statistical exhaustiveness, Appl. Math. Lett. \textbf{25} 10 (2012) 1447-1451.
\bibitem{CasertaMaioKocinacStatisticalConvergenceinFunctionSpaces} A. Caserta, G. Di Maio, and Ljubisa D.R. Kocinac, Statistical Convergence in Function Spaces, Abstr. Appl. Anal. \textbf{2011}, (2011) Article ID 420419 11 pages doi:10.1155/2011/420419
\bibitem{ConnorandErdmannSequentialdefinitionsofcontinuityforrealfunctions}  J. Connor, and K.-G. Grosse-Erdmann, Sequential definitions of continuity for real functions, Rocky Mountain J. Math. \textbf{33} 1 (2003) 93-121.
\bibitem{CakalliLacunarystatisticalconvergenceintopgroups}  H. \c{C}akall\i,  Lacunary statistical convergence in topological groups,
  Indian J. Pure Appl. Math. \textbf{26} 2 (1995) 113-119.
\bibitem{CakalliSlowlyoscillatingcontinuity} H. \c{C}akall\i,  Slowly oscillating continuity, Abstr. Appl. Anal. \textbf{2008}, Article ID 485706 5 pages (2008) doi:10.1155/2008/485706
\bibitem{CakalliSequentialdefinitionsofcompactness} H. \c{C}akall\i,  Sequential definitions of compactness, Appl. Math. Lett. \textbf{21} 6 (2008) 594-598.
\bibitem{CakalliAstudyonstatisticalconvergence}  H. \c{C}akall\i,  A study on statistical convergence, Funct. Anal. Approx. Comput. \textbf{1} 2 (2009) 19-24.
\bibitem{CakalliForwardcontinuity} H. \c{C}akall\i,  Forward continuity,  J. Comput. Anal. Appl. \textbf{13} 2 (2011) 225-230.
\bibitem{CakalliOnGcontinuity} H. \c{C}akall\i,  On $G$-continuity, Comput. Math. Appl. \textbf{61} 2 (2011) 313-318.
\bibitem{CakalliDeltaquasiCauchysequences} H. \c{C}akall\i,  $\delta$-quasi-Cauchy sequences, Math. Comput. Modelling. \textbf{53} 1-2 (2011) 397-401.
\bibitem{CakalliOnDeltaquasislowlyoscillatingsequences}  H. \c{C}akall\i,  On $\Delta$-quasi-slowly oscillating sequences, Comput. Math. Appl. \textbf{62} 9  (2011) 3567-3574.
\bibitem{CakalliStatisticalquasiCauchysequences}  H. \c{C}akall\i,  Statistical quasi-Cauchy sequences, Math. Comput. Modelling. \textbf{54} 5-6 (2011) 1620-1624. \bibitem{CakalliStatisticalwardcontinuity} H. \c{C}akall\i,  Statistical ward continuity, Appl. Math. Lett. \textbf{24} 10 (2011) 1724-1728.
\bibitem{CakalliSequentialdefinitionsofconnectedness} H. \c{C}akall\i, Sequential definitions of connectedness, Appl. Math. Lett. \textbf{25} 3 (2012), 461-465.
\bibitem{CakalliAvariationonwardcontinuity} H. \c{C}akall\i,  A variation on ward continuity, Filomat. \textbf{27} 8 (2013) 1545-1549.
\bibitem{CakalliVariationsonquasiCauchysequences} H. \c{C}akall\i,  Variations on quasi-Cauchy sequences, Filomat, \textbf{29} 1 (2015) 13-19.
\bibitem{CakalliUpwardanddownwardstatisticalcontinuities} H. \c{C}akall\i,  Upward and downward statistical continuities, Filomat, \textbf{29} 10 (2015) 2265-2273.
\bibitem{CakalliandPratulFuzzycompactnessviasummability}  H. \c{C}akall\i,  and  P. Das, Fuzzy compactness via summability, Appl. Math. Lett. \textbf{22} 11  (2009) 1665-1669.
\bibitem{CakalliandErsanStronglylacunarywardcontinuityin2normedspaces} H. \c{C}akall\i\;, S. Ersan, Strongly Lacunary Ward Continuity in 2-Normed Spaces, Scientific World Journal, \textbf{2014}, Article ID 479679, (2014) 5 pages.  doi: http://dx.doi.org/10.1155/2014/479679
\bibitem{CakalliandErsanLacunarywardcontinuityin2normedspaces} H. Cakalli and S. Ersan, Lacunary Ward Continuity in 2-normed Spaces, Filomat, \textbf{29} 10 (2015), 2257-2263.
\bibitem{CakalliandKaplanAstudyonNthetaquasiCauchysequences} H. \c{C}akall\i,  and H. Kaplan, A study on N-theta quasi-Cauchy sequences,  Abstr. Appl. Anal. \textbf{2013}, Article ID 836970 4 pages  (2013)
\bibitem{CakalliandKhan} H. \c{C}akall\i,  and M.K. Khan,  Summability in Topological Spaces, Appl. Math. Lett. \textbf{24}, (2011) 348-352.
\bibitem{CakalliandMucukconnectednessviaasequentialmethod} H. Cakalli and O. Mucuk, On connectedness via a sequential method, Revista de la Union Matematica Argentina \textbf{54} 2 (2013), 101-109.
\bibitem{CakalliandSavasStatisticalconvergenceofdoublesequences} H. \c{C}akalli  and E. Sava\c{s}, Statistical convergence of double sequences in topological groups, J. Comput. Anal. Appl. \textbf{12} 2 (2010) 421-426.
\bibitem{CakalliandSonmezSlowlyoscillatingcontinuityinabstractmetricspaces} H. \c{C}akalli and A. Sonmez, Slowly oscillating continuity in abstract metric spaces, Filomat \textbf{27} (2013) 925-930.
\bibitem{CakalliandSonmezandArasLamdastatisticallywardcontinuity} H. \c{C}akalli, A. Sonmez, and C.G. Aras, $\lambda$-statistically ward continuity, An. Stiint. Univ. Al. I. Cuza Iasi. Mat. (N.S.) Tomul LXIII, f. 2 (2017) 313-321. http://www.math.uaic.ro/~annalsmath/pdf-uri
\bibitem{CakalliandSonmezandGenc}  H. \c{C}akall\i, A. Sonmez, and  C. Genc, On an equivalence of topological vector space valued cone metric spaces and metric spaces, Appl. Math. Lett. \textbf{25} (2012), 429-433.
\bibitem{CanakandDikNewtypesofcontinuities}  I. Canak, and M. Dik, New types of continuities, Abstr. Appl. Anal. \textbf{2010}, Article ID 258980 (2010) 6 pages.  \bibitem{FDikMDikandCanak}  F. Dik,  M. Dik and I. Canak,  Applications of subsequential Tauberian theory to classical Tauberian theory, Appl. Math. Lett. \textbf{20} 8 (2007) 946-950.
\bibitem{DjurcicandKocinacandZizovicDoublesequencesandselections} D. Djurcic, Ljubisa D.R. Kocinac, M.R. Zizovic,  Double sequences and selections, Abstr. Appl. Anal.  \textbf{2012}, Art. ID 497594 (2012) 6 pages.
\bibitem{ErsanandCakalliWardContinuityin2NormedSpaces} S. Ersan and H. Cakalli, Ward Continuity in 2-Normed Spaces, Filomat \textbf{29} 7 (2015), 1507-1513. DOI 10.2298/FIL1507507E
\bibitem{Fast}  Fast, H.  Sur la convergence statistique, Colloq. Math. \textbf{2} (1951) 241-244.
\bibitem{Fridy} Fridy, J.A.   On statistical convergence. Analysis, \textbf{5} (1985) 301-313.
\bibitem{FridyandKhan} Fridy, J.A. and Khan, M.K.  Statistical extensions of some classical Tauberian theorems, Proc. Amer. Math. Soc.  \textbf{128} 8 (2000) 2347-2355.
\bibitem{FridyandOrhanlacunarystatisconvergence} J.A. Fridy, and C. Orhan, Lacunary statistical convergence, Pacific J. Math. \textbf{160} 1 (1993) 43-51.
\bibitem{KocinacSelectionpropertiesinfuzzymetricspaces}  Ljubisa D.R. Kocinac,  Selection properties in fuzzy metric spaces, Filomat. \textbf{26} 2 (2012) 305-312.
\bibitem{MucukandCakalliGsequentiallyconnectednessfortopologicalgroupswithoperations} O. Mucuk and H. Cakalli, {\em G-sequentially connectedness  for topological groups with operations}, AIP Conference Proceedings, \textbf{1759}, 020038, (2016). doi: 10.1063/1.4959652
\bibitem{MucukandSahanOnGsequentialcontinuity} O. Mucuk, T. \c{S}ahan¸ On $G$-Sequential Continuity, Filomat \textbf{28} 6 (2014), 1181-1189.
\bibitem{MursaleenandAlotaibiOnIconvergenceinrandom2normedspaces} M. Mursaleen, A. Alotaibi, On I-convergence in random 2-normed spaces, Math. Slovaca \textbf{61} 6 (2011), 933-940
\bibitem{PalandSavasandCakalliIconvergenceonconemetricspaces} S.K. Pal,  E. Savas, and H. Cakalli, $I$-convergence on cone metric spaces, Sarajevo J. Math. \textbf{9} (2013), 85-93.
\bibitem{PattersonandCakalliQuasiCauchydoublesequences} R.F. Patterson and H. Cakalli, Quasi Cauchy double sequences, Tbilisi Mathematical Journal \textbf{8} 2 (2015), 211-219.
\bibitem{SavasOnlacunarystrongsigmaconvergence} E. Sava\c{s}, A study on absolute summability factors for a triangular matrix, Math. Inequal. Appl. \textbf{12} 1 (2009) 141-146.
\bibitem{SavasandNurayOnsigmastatisticallyconvergenceanlacunarysigmastatisticallyconvergence}   E. Sava\c{s} and F. Nuray, On $\sigma$-statistically convergence and lacunary $\sigma$-statistically convergence, Math. Slovaca. \textbf{43} 3 (1993) 309-315.

\bibitem{PowellandShahSummabilityTheoryandItsApplications} R.E. Powell, S.M. Shah, Summability Theory and Its Applications, Prentice-Hall of India, New Delhi (1988)
\bibitem{SonmezOnparacompactnessinconemetricspaces} A. Sonmez, {\em On paracompactness in cone metric spaces}, Appl. Math. Lett. \textbf{23}, 494-497, (2010).
\bibitem{SonmezandCakalliConenormedspacesandweightedmeans} A. Sonmez, and H. Cakalli, Cone normed spaces and weighted means, Math. Comput. Modelling  \textbf{52}  (2010), 1660-16660.
\bibitem{UnverAbelsummabilityintopologicalspaces} M. \"Unver, Abel summability in topological spaces, Monatsh Math 178 (2015) 633-643. https://doi.org/10.1007/s00605-014-0717-0
\bibitem{Vallin} Vallin, R.W.  Creating slowly oscillating sequences and slowly oscillating continuous functions With an appendix by Vallin and H. Cakalli, Acta Math. Univ. Comenianae. \textbf{25} 1 (2011) 71-78.
\bibitem{YildizIstatistikselboslukludelta2quasiCauchydizileri} \c{S}. Y\i ld\i z, {\em \.{I}statistiksel bo\c{s}luklu delta 2 quasi Cauchy dizileri}, Sakarya University Journal of Science, \textbf{21}, 6, (2017). DOI: 10.16984/saufenbilder.336128 , http://www.saujs.sakarya.edu.tr/issue/26999/336128










\end{thebibliography}
\end{document}